\documentclass[11pt]{article}
\usepackage{graphicx}
\usepackage{amsmath}
\usepackage{amsfonts}
\usepackage{epsfig}
\usepackage{setspace}

\newcommand{\be}{\begin{equation}}
\newcommand{\ee}{\end{equation}}
\newcommand{\beqn}{\begin{eqnarray}}
\newcommand{\eeqn}{\end{eqnarray}}
\newcommand{\beqns}{\begin{eqnarray*}}
\newcommand{\eeqns}{\end{eqnarray*}}

\newcommand{\card}{\mbox{card}}
\newcommand{\Var}{\mbox{Var}}

\newcommand{\EE}{\ensuremath{{\mathbb E}}}
\newcommand{\II}{\ensuremath{{\mathbb I}}}

\newcommand{\fr}[1]{(\ref{#1})}

\newcommand{\Te}{\Theta}

\newtheorem{lemma}{Lemma}
\newtheorem{theorem}{Theorem}

\newtheorem{remark}{Remark}

\begin{document}

\title{\Large{\bf Estimation in nonparametric regression model with additive and multiplicative noise via Laguerre series}}

\author{
\large{ Rida Benhaddou}  \footnote{E-mail address: Benhaddo@ohio.edu}
  \\ \\
Department of Mathematics, Ohio University, Athens, OH 45701} 
\date{}

\doublespacing
\maketitle
\begin{abstract}
We look into the nonparametric regression estimation with additive and multiplicative noise and construct adaptive thresholding estimators based on Laguerre series. The proposed approach achieves asymptotically near-optimal convergence rates when the unknown function belongs to Laguerre-Sobolev space. We consider the problem under two noise structures; (1) { i.i.d.} Gaussian errors and (2) long-memory Gaussian errors. In the { i.i.d.} case, our convergence rates are similar to those found in the literature. In the long-memory case, the convergence rates depend on the long-memory parameters only when long-memory is strong enough in either noise source, otherwise, the rates are identical to those under { i.i.d.} noise. \\


{\bf AMS (2000) Subject Classification: 62G05, 62G20, 62G08 }
 \end{abstract} 

\section{Introduction}

Consider a nonparametric regression model with both multiplicative and additive noise  
\be
y_i=f(t_i)\varepsilon_i +\sigma z_i, \ \ \ i=1, 2, \cdots, N, \label{conveq}
\ee
where $\varepsilon_i$ and $z_i$ are zero-mean $(1)$ independent and identically distributed Gaussian random variables with variance equal to 1 and (2) $\varepsilon_i$ and $z_i$ are Gaussian with long-memory structure, and $\sigma$ is a known  positive constant. The function $f(t)$ is the unknown response, it is real-valued and is defined on the interval $[0, b]$ with $b>0$ a fixed real number.  In addition, $t_i$,  $i=1, 2, \cdots, N$, are independent and identically distributed random variables, drawn from a known probability density function $g$ with support $[0, b]$. It is assumed that the quantities $t_i$, $\varepsilon_i$ and $z_i$ are independent from one another for any  $i\in \{1, 2, \cdots, N\}$. The goal is to estimate $h(t)=f^2(t)$ based on data points $(t_1, y_1), (t_2, y_2), \cdots, (t_N, y_N)$. 

This problem, under various settings, has been studied considerably by the means of a number of nonparametric methods, including kernel smoothing, splines and wavelets, and the list of articles includes, in chronological order, Hardle and Tsybakov~(1997), Brown and Levine~(2007), Cai and Wang~(2008), Kulik and Wichelhaus~(2011) and Chichignoud~(2012). Most recently, Chesneau, El Kolei, Kou and Navarro~(2020) considered the problem in a multivariate setting  and proposed a wavelet thresholding approach to solve it. This problem has a great deal of applications, for instance, { in Global Positioning System (GPS) signal propagation modeling where there is empirical evidence that in heavy multi-path urban areas, the GPS signal encounters both additive and multiplicative noise (see Huang et al.~(2013)), or} in finance where one is interested in estimating the variance from the returns of an asset and the interested reader may refer to Chesneau et al.~(2020) for more.  Almost all of these articles  assume that the error terms are white noise processes or {\bf i.i.d.} noise. However, empirical evidence has shown that, even at large lags, the correlation structure in the errors may take the power-like form. This phenomenon is referred to as long-memory (LM) or long-range dependence (LRD).

Long-memory has been investigated quite considerably in  many nonparametric estimation problems, including regression and deconvolution and the list includes Wang~(1996, 1997), Comte, Dedecker and Taupin~(2008),  Kulik and Raimondo, M. (2009), Kulik and Wichelhaus~(2011), Wishart~(2013), Benhaddou, Kulik, Pensky and Sapatinas~(2014), Benhaddou~(2016), Benhaddou~(2018a, 2018b) and  Benhaddou and Liu~(2019).  

The application of Laguerre series to Nonparametric estimation has become popular as of late and the list includes the application to density estimation in Comte et al.~(2008) and Comte and Genon-Catalot~(2015), the estimation of linear functionals of a density function in  Mabon~(2016) and Laplace deconvolution in Vareschi~(2015), Comte, Cuenod, Pensky, and Rozenholc~(2017) and Benhaddou, Pensky and Rajapakshage~(2019). 

The objective of the paper is to solve the nonparametric regression model with both multiplicative and additive i.i.d., and long-memory Gaussian noise via Laguerre hard-thresholding when the design points are random and follow known probability density function $g$. We derive lower bounds for the $L^2$-risk when $h=f^2$ belongs to some Laguerre-Sobolev ball of radius $A >0$, and then construct an adaptive Laguerre-thresholding estimator  for $h=f^2$. In addition, we show that the proposed estimator attains asymptotically near-optimal convergence rates. Furthermore, we demonstrate that long-memory has a detrimental effect on the convergence rates only when it is strong enough in either noise source. In which case, the convergence rates depend on the smoothness of the unknown function $h=f^2$ and the long-memory parameter associated with the stronger dependence between the two noise sources. It turns out that the present rates are identical to those in Chesneau et al.~(2020)  with $d=1$ for the i.i.d. case and with $\alpha_1=\alpha_2=1$ for the long-memory case.  Similarly, our rates are comparable to those in Brown and Levine~(2007) and Cai and Wang~(2008)  in their treatment of the regression variance estimation when the unknown mean function is  smooth enough. 

 \section{Estimation Algorithm}
 For the rest of the paper, let $\|h\|$ denote the $L^2([0, \infty))$-norm of the function $h$. Given a matrix $A$, let $A^T$ be its transpose, $\lambda_{\max}(A)$ be its largest eigenvalue in magnitude, $\|A\|_F=\sqrt{Tr\left(A^TA\right)}$ and $\|A\|_{sp}=\lambda_{\max}\left(A^TA\right)$ be, respectively, its Frobenius and the spectral norms. In addition, let $(a\vee b)=\max(a, b)$ and $(a\wedge b)=\min(a, b)$. Consider the orthonormal basis that consists of the system of Laguerre functions 
\be
\varphi_k(t)= e^{-t/2}L_k(t), \ \ k=0, 1, \cdots, 
\ee
where $L_k(t)$ are Laguerre polynomials (see. e.g., Gradshtein and Ryzhik~(1980), Section 8.97). Since the functions $\varphi_k(t)$, $k=0, 1, \cdots,$ form an orthonormal basis on $[0, \infty)$, the function $h(.)=f^2(.)$ can be expanded over this basis as follows
\be \label{lagmattr}
h(t)\equiv f^2(t)=\sum^{\infty}_{k=0}\theta_k\varphi_k(t),
\ee \label{theta-h}
where  $\theta_k=\int^{b}_0h(t)\varphi_k(t)dt=\int^{b}_0f^2(t)\varphi_k(t)dt$. Under i.i.d noise case, similar to  Chesneau et al.~(2020), an estimator for $\theta_k$ is given by
\be \label{theta-h}
\widehat{\theta}_k=\frac{1}{N}\sum^N_{i=1}\left[y_i^2\frac{\varphi_k(t_i)}{g(t_i)}-\sigma^2\int^{b}_0\varphi_k(t)dt\right]\II\left(\Omega_k(i) \right),
\ee
where $\Omega_k(i)=\left\{i: \left|y_i^2\frac{\varphi_k(t_i)}{g(t_i)}-\sigma^2\int^{b}_0\varphi_k(t)dt\right|  \leq   \sqrt{\frac{N}{\ln(N)}} \right \}$. Similarly, under long-memory noise case, an unbiased estimator for $\theta_k$ is given by
\be \label{theta-h-lm}
\widehat{\theta}_k=\frac{1}{N}\sum^N_{i=1}\left[y_i^2\frac{\varphi_k(t_i)}{g(t_i)}-\sigma^2\int^{b}_0\varphi_k(t)dt\right].
\ee
Then, consider the hard-thresholding estimator for $h=f^2$
\be\label{est-f}
\widehat{h}_M(t)= \sum^{M-1}_{l=0}\widehat{\theta}_l\II\left(|\widehat{\theta}_l| > \lambda_l\right)\varphi_l(t),
\ee 
where the quantities $M$ and $\lambda_l$ will be determined under the two different setups in the proceeding sections.\\
Next is the list of conditions that will be utilized in the derivation of the theoretical results. \\
{\bf Assumption A.1.} $f \in {\bf{L^{2}}}\left[0, { b} \right)$ is bounded above, that is, there exists positive constant $M_2<\infty$ such that  $f(t) \leq M_2$, for all  $t \in \left[0, { b} \right)$.\\
{\bf Assumption A.2.} The probability density function $g$ is uniformly bounded, that is, on  $\left[0, { b} \right)$ there exist positive constants $m_1$ and $m_2$, with $0< m_1\leq m_2 < \infty$, such that $m_1 \leq g(t) \leq m_2$. 
 \begin{remark}	  
{ Assumption {\bf A.2.} is valid for instance when  $g$ is the uniform distribution. In such case $m_1=m_2=1/b$. If $g$ is not bounded below, such as in the case of beta distribution with b=1, a variation of the present procedure will be needed and this would be another direction for future research. The idea is to consider the generalized Laguerre function basis instead, which is defined by 
\be
\varphi^{(a)}_k(x)={\left[\frac{k!}{\Gamma(k+a+1)}\right]^{1/2}} e^{-x/2}x^{a/2}L^{(a)}_k(x), \ \ k=0, 1, \cdots, 
\ee
where $L^{(a)}_k(t)$ are generalized Laguerre polynomials with parameter $a$, $a> 0$ (see. e.g., Gradshtein and Ryzhik~(1980), Section 8.97), and select the parameter $a$ according to $g$ at hand}.
  \end{remark}
{\bf Assumption A.3.} \label{A7} The function $h(t)=f^2(t)$ belongs to a Laguerre-Sobolev space. In particular, Laguerre coefficients of $h$, $\theta_{l}$ satisfy
\be  \label{eq11}
 {\bf B}^{s}(A)=\left \{ h \in L^2[0, b]: \sum^{\infty}_{l=0} (l\wedge1)^{2s} \theta^2_l \leq A\right \}.
\ee
 we are in the position to fill in the details of the estimator and find the minimax lower bound for the quadratic risk and compare it to asymptotic upper bound for the mean squared error of our estimator. 
  \begin{remark}	  
{ Functional spaces of type \fr{eq11} have been introduced in Bongioanni and Torrea~(2009) to study Laguerre operators, and the connection with Laguerre coefficients was established in Comte and Genon-Catalot~(2015)}.  
  \end{remark}
\section{Asymptotic minimax and adaptivity: the i.i.d. case}
  \noindent
We define the minimax $L^2$-risk over a set $\Theta$ as 
$$R(\Theta)=\inf_{\tilde{h}}\sup_{h\in \Theta}\EE\|\tilde{h}-h\|^2,$$
where the infimum is taken over all possible estimators $\tilde{h}$ of $h$. 
\begin{theorem}\label{th:lowerbds} Let  Assumptions {\bf{A.1-A.3}} hold. Then, as $N\rightarrow \infty$,
 \be \label{lowerbds}
 R({\bf{B^s}}(A))\geq CA^2
 \left[ \frac{1}{A^2N}\right]^{\frac{2s}{2s +1}}.
   \ee
 \end{theorem}
 \begin{lemma} \label{lem:Var}
Let conditions ${\bf{A.1}}$ and ${\bf{A.2}}$ hold and let $\widehat{\theta}_l$ be defined in \fr{theta-h}. Then, for $l= 1, 2, \cdots, M-1$, as $N \rightarrow \infty$, one has 
\be \label{var-bias}
\EE|\widehat{\theta}_{l}-\theta_l|^2 = O\left(\frac{1}{N}\right).
\ee
\end{lemma}
Based on Lemma 1, choose the thresholds $\lambda_l$ such that 
\be \label{thresh1}
\lambda_l=\gamma\frac{\sqrt{\ln(N)}}{\sqrt{N}}.
\ee
 In addition choose the truncation level $M$ as
\be  \label{Lev:Jad}
M=N.
\ee
 \begin{lemma} \label{lem:Lardev}
Let conditions ${\bf{A.1}}$ and ${\bf{A.2}}$ hold and let $\widehat{\theta}_l$ be defined in \fr{theta-h}. Then, for $l= 1, 2, \cdots, M-1$, if $\rho\gamma >1$, as $N \rightarrow \infty$, one has 
\be  \label{Largdev}
\Pr \left(| \widehat{\theta}_{l}-\theta_l |> \rho\lambda_l\right)= O \left ( \left[\frac{1}{N}\right]^{\tau}  \right),
\ee
where $\tau$ is a positive parameter that is large enough and $\rho$ is such that $0 < \rho <1$. 
\end{lemma}
\begin{theorem} \label{th:upperbds-2}
Let $s \geq 1/2$ and let  $\tilde{h}_{M}(t)$ be the Laguerre estimator defined in \fr{est-f} with $M$ given in \fr{Lev:Jad} and $\lambda_l$ given in \fr{thresh1}. Suppose assumptions ${\bf{A.1}}$-${\bf{A.3}}$ hold. Then,  if $\tau$ is large enough, as $N \rightarrow \infty$, one has
 \be \label{upperbds-2}
 \sup_{h\in B^s(A)} \EE\|\widehat{h}_M-h\|^2\leq CA^2
\left[\frac{ \ln(N)}{A^2N}\right]^{\frac{2s}{2s+1}}. 
 \ee
\end{theorem}
 \begin{remark}	
 Theorems \ref{th:lowerbds} and \ref{th:upperbds-2} imply that, for the $L^2$-risk, the estimator \fr{est-f} with $\lambda_l$ given by \fr{thresh1} and $M$ chosen according to \fr{Lev:Jad} is adaptive and asymptotically near-optimal, within a logarithmic factor of $N$, over all  Laguerre-Sobolev balls $B^{s}(A)$. 
  \end{remark} 
  \begin{remark}	  
The convergence rates match those in Chesneau et al.~(2020), with $d=1$ in their treatment of the problem using wavelets when the function under consideration belongs to a certain Besov ball. 
  \end{remark}
   \begin{remark}	
 Our rates are comparable to those in Brown and Levine~(2007) and Cai and Wang~(2008)  in their treatment of the regression variance estimation when the unknown mean function is  smooth enough.
  \end{remark}
  \begin{remark}	
 In GPS signal detection application, the size of $\sigma$  in equation \fr{conveq} will dictate whether the multiplicative noise will be considered or ignored in the analysis (see Huang et al.~(2013)). In addition, $\sigma$ may not be know in advance but it can be estimated from the data. Providing fully data-driven procedure is beyond the scope of this work so we assume $\sigma$ is known. 
  \end{remark}
   \section{Asymptotic minimax and adaptivity: the long-memory case}
Let $\pmb{ \varepsilon}_N$ be the random vector with elements $\varepsilon_i$, $i=1, 2, \cdots, N$, and covariance matrix $\Sigma_1=Cov(\pmb{\varepsilon}_N)=\EE\left[\pmb{\varepsilon}_N \pmb{\varepsilon}^T_N\right]$, and let $\pmb{ z}_N$ be the random vector with elements $z_i$, $i=1, 2, \cdots, N$ and covariance matrix $\Sigma_2=Cov(\pmb{z}_N)=\EE\left[\pmb{z}_N \pmb{z}^T_N\right]$. \\
 {\bf Assumption A.4.} The vectors $\pmb{\varepsilon}_N$ and $\pmb{ z}_N$ allow the decomposition 
\be \label{epseta}
\pmb{\varepsilon}_N=A_1\pmb{\eta}^{(1)}_N, \ \ \ \pmb{ z}_N=A_2\pmb{\eta}^{(2)}_N,
\ee
where $\pmb{\eta}^{(j)}_N$, $j=1, 2$, is a random vector with zero-mean independent Gaussian $\eta^{(j)}_i$ having equal variance, $i=1, 2, \cdots, N$  and $A_j$ is some non-random matrix. $\{\varepsilon_{i}\}_{i\geq 1}$ and $\{z_i\}_{i\geq 1}$ are zero-mean,  stationary Gaussian sequences with auto-covariance functions $\gamma_1(h)=\EE\left[\varepsilon_i\varepsilon_{i+h}\right]$ and $\gamma_2(h)=\EE\left[z_iz_{i+h}\right]$, satisfying 
\be \label{auto-cov}
\gamma_1(h)\asymp h^{-\alpha_1},\ \  \gamma_2(h)\asymp h^{-\alpha_2}.
\ee
\noindent
{\bf Assumption A.5.}  There exist constants $c^{(j)}_1$ and $c^{(j)}_2$, $j=1, 2$, ($0 < c^{(j)}_1 \leq c^{(j)}_2 < \infty$), independent of $N$, such that
\be
c^{(j)}_1 n^{1-\alpha_j}\leq \lambda_{min}(\Sigma_j) \leq \lambda_{max}(\Sigma_j)\leq c^{(j)}_2 n^{1-\alpha_j}, \ \ 0< \alpha_j \leq 1, \label{LRD}
\ee
where $\alpha_j$, $j=1, 2$, are the long-memory parameters associated with the matrices $\Sigma_j$, respectively.
\begin{theorem}\label{th:lowerbds-lm} Let  Assumptions {\bf{A.1-A.3}} and {\bf{A.5}} hold. Then, provided that $f$ is bounded away from zero, as $N\rightarrow \infty$, one has
 \be \label{lowerbds-lm}
 R({\bf{B^s}}(A))\geq CA^2
\left\{ \begin{array}{ll}  
\left[\frac{1}{A^2N}\right]^{\frac{2s}{2s+1}}, \ & \mbox{if}\ (\alpha_1\wedge \alpha_2) \geq 1/2, \\
  \left[\frac{ 1}{A^2(N^{2\alpha_1}\vee N^{2\alpha_2})}\right]^{\frac{2s}{2s+1}} & \mbox{if}\ otherwise. 
    \end{array} \right.
   \ee
 \end{theorem}
  \begin{lemma} \label{lem:Var-lm}
Let conditions ${\bf{A.1}}$, ${\bf{A.2}}$, ${\bf{A.4}}$ and ${\bf{A.5}}$ hold and let $\widehat{\theta}_l$ be defined in \fr{theta-h-lm}. Then, for $l= 1, 2, \cdots, M-1$, as $N \rightarrow \infty$, one has 
\be \label{var-bias2}
\EE|\widehat{\theta}_{l}-\theta_l|^2 = \left\{ \begin{array}{ll} 
O\left(\frac{1}{N}\right)\ & \mbox{if}\ (\alpha_1\wedge \alpha_2) \geq 1/2, \\
O\left(\frac{1}{N^{2\alpha_1}}\vee \frac{1}{N^{2\alpha_2}}\right)\ & \mbox{if}\ otherwise.
  \end{array} \right.
\ee
\end{lemma}
Based on Lemma 3, choose the thresholds $\lambda_l$ such that 
\be \label{thresh2}
\lambda_l=\left\{ \begin{array}{ll}  
\gamma\frac{\sqrt{\ln(N)}}{\sqrt{N}}, \ & \mbox{if}\ (\alpha_1\wedge \alpha_2) \geq 1/2, \\
\gamma_1\frac{\ln{(N)}}{N^{\alpha_1}}\vee \gamma_2 \frac{\ln{(N)}}{N^{\alpha_2}}\ & \mbox{if}\ otherwise.
  \end{array} \right.
  \ee
 In addition, choose the maximal level $M$ 
\be  \label{Lev:Jad2}
M=\left\{ \begin{array}{ll}  
N, \ & \mbox{if}\ (\alpha_1\wedge \alpha_2) \geq 1/2, \\
{N^{2\alpha_1}}\wedge{N^{2\alpha_2}}\ & \mbox{if}\ otherwise.
  \end{array} \right.
\ee
 \begin{lemma} \label{lem:Lardev-lm}
Let conditions ${\bf{A.1}}$, ${\bf{A.2}}$, ${\bf{A.4}}$ and ${\bf{A.5}}$ hold and let $\widehat{\theta}_l$ be defined in \fr{theta-h-lm}. Then, for $l= 1, 2, \cdots, M-1$, as $N \rightarrow \infty$, one has 
\be  \label{Largdev}
\Pr \left(| \widehat{\theta}_{l}-\theta_l |> \rho\lambda_l\right)= O \left ( \left[\frac{1}{N}\right]^{\tau}  \right),
\ee
where $\tau$ is a positive parameter that is large enough and $\rho$ is such that $0 < \rho <1$. 
\end{lemma}
\begin{theorem} \label{th:upperbds-lm}
Let $s \geq 1/2$ and let  $\tilde{h}_{M}(t)$ be the Laguerre estimator defined in \fr{est-f} with $M$ given in \fr{Lev:Jad2} and $\lambda_l$ given in \fr{thresh2}. Suppose assumptions ${\bf{A.1}}$-${\bf{A.5}}$ hold. Then,  if $\tau$ is large enough, as $N \rightarrow \infty$, one has
 \be \label{upperbds-lm}
 \sup_{h\in B^s(A)} \EE\|\widehat{h}_M-h\|^2\leq CA^2\left\{ \begin{array}{ll}  
\left[\frac{ \ln(N)}{A^2N}\right]^{\frac{2s}{2s+1}}, \ & \mbox{if}\ (\alpha_1\wedge \alpha_2) \geq 1/2, \\
  \left[\frac{ \ln^2(N)}{A^2(N^{2\alpha_1}\vee N^{2\alpha_2})}\right]^{\frac{2s}{2s+1}} & \mbox{if}\ otherwise. 
    \end{array} \right.
 \ee
\end{theorem}
\begin{remark}	
Notice that when both long-memory parameters are large enough, in particular when $(\alpha_1 \wedge \alpha_2) \geq 1/2$, the convergence rates are identical to those under i.i.d. errors. In such case they match those in  Chesneau et al.~(2020) directly if, in their notation, $d=1$.
 \end{remark}
\begin{remark}	
When the long-memory is strong, which corresponds to relatively low $\alpha_j$, the convergence rates depend on the smoothness of the unknown function $h=f^2$ and the long-memory parameter associated with the stronger dependence, $(\alpha_1 \wedge \alpha_2)$,  between the two noise sources. These rates are completely new and provide an extension of the problem in a different direction. 
  \end{remark}
\begin{remark}	
Notice that the i.i.d. case can also be handled using estimators \fr{theta-h-lm} along with the choice of thresholds \fr{thresh1} and truncation level $M$ based on \fr{Lev:Jad} and achieve the same convergence rates. 
  \end{remark}
 \section{Proofs }
In order to prove Theorem \ref{th:lowerbds}, we use the following lemma  
\begin{lemma} (Lemma $A.1$ of Bunea et al.~(2007)) 
\label{lem:Bunea} 
Let $\Te$ be a set of functions of cardinality $\card(\Te)\geq 2$ such that\\
(i) $\|f-g\|_p^p \geq 4\delta^p, \ for\  f, g \in \Te, \ f \neq g, $\\
(ii) the Kullback divergences $K(P_f, P_g)$ between the measures $P_f$ and $P_g$ 
satisfy the inequality $K(P_f, P_g) \leq \log(\card(\Te))/16,\ for\ f,\ g \in \Te$.\\
Then, for some absolute positive constant $C_1$, one has 
$$
 \inf_{f_n}\sup_{f\in \Te} \EE_f\|f_n-f\|_p^p \geq C_1 \delta^p,
$$
where $\inf_{f_n}$ denotes the infimum over all estimators.
\end{lemma}
 {\bf Proof of Theorem \ref{th:lowerbds} and Theorem \ref{th:lowerbds-lm}}. Let $\omega$ be the vector with components $\omega_l \in \{ -1, 1\}$, $l=0, 1, \cdots, L-1$, and denote the set of all possible values of $\omega$ by $\Omega$. Let $h_L$ be the functions of the form
\be  \label{f_j}
h_{L}(t)=\rho_L\sum^{L-1}_{l=1}\omega_l\varphi_l(t),\ \ \omega_l \in \{ -1, 1\}.
\ee
Observe that $\omega$ has $L$ and therefore $\Omega$ will have cardinality $\card(\Omega)=2^L$. By \fr{eq11}, it is easy to verify that  $h_{L}(t) \in B^{s}(A)$ with the choice $\rho^2_L=C A^2L^{-(2s+1)}$. Take $\tilde{h}_{L}$ of the form of \fr{f_j} but with $\tilde{\omega}_l\in \{ -1, 1\}$, then applying Varshamov-Gilbert Lemma (\cite{tsybakov}, p 104), the $L^2$-norm of the difference is
\be
\| h_{L}(t)-\tilde{h}_{L}(t)\|^2 \geq  \frac{L\rho^2_L}{8}.
\ee
To prove {\bf Theorem \ref{th:lowerbds}}, define the quantities $\pmb{h}_i=h_L (t_i)+ f_L(t_i) \nu_i+ \sigma_1\mu_i$, $i=1, 2, \cdots, N$, where $\{\nu_{i}\}_{i\geq 1}$ and $\{\mu_i\}_{i\geq 1}$ are i.i.d. $N(0, 1)$ sequences that are independent of each other. Let $P_{h_L}$ be the probability law of the process $\pmb{h}_i$ under the hypothesis $h_L$ defined in \fr{f_j}. Then, by Assumptions ${\bf A.2}$ and ${\bf A.3}$, the Kullback divergence can be written as 
\beqn
K(P_{h_{L}}, P_{\tilde{h}_{L}})&\leq& {N\rho^2_L}\sum^{L-1}_{l=0}|\tilde{\omega}_l-\omega_l|^2\frac{\int^{b}_0\varphi^2_l(t)g(t)dt}{2(\int^{b}_0\tilde{f}^2_L(t)g(t)dt+\sigma_1^2)}\nonumber\\
& \leq & \frac{N\rho^2_L}{\sigma_1^2} 2L \max_{l \leq L-1}|\varphi_l(t)|^2=CNA^2L^{-(2s+1)}L.
\eeqn
Now, to apply {\bf Lemma \ref{lem:Bunea}}, choose 
\be \label{klblam}
A^2 L^{-(2s+1)}LN\leq \pi_0L.
\ee
Hence, the proof is complete by taking  
\be \label{Lop}
L = C\left[{A^2N}\right]^{\frac{1}{2s +1}}. \ \Box
\ee
To prove {\bf Theorem \ref{th:lowerbds-lm}}, define the vectors $\pmb{h}_N$ whose elements are quantities $\pmb{h}_i=h_L (t_i)+ f_L(t_i) \nu_i+ \sigma_1\mu_i$, $i=1, 2, \cdots, N$, such that $\pmb{h}_N\sim N({\bf h}_L, \tilde{\Sigma}_1+\sigma^2_1\tilde{\Sigma}_2)$. Here, $\{\nu_{i}\}_{i\geq 1}$ and $\{\mu_i\}_{i\geq 1}$ are zero-mean,  stationary Gaussian sequences with auto-covariance functions $\gamma_1(h)=\EE\left[\nu_i\nu_{i+h}\right]$ and $\gamma_2(h)=\EE\left[\mu_i\mu_{i+h}\right]$, satisfying 
\be \label{auto-cov2}
\gamma_1(h)\asymp h^{-2\alpha_1},\ \  \gamma_2(h)\asymp h^{-2\alpha_2}.
\ee
Notice that under \fr{auto-cov2}, (1) $\lambda_{min}(\tilde{\Sigma}_i) \asymp 1$, $i=1, 2$,  if $(\alpha_1\wedge \alpha_2) \geq 1/2$, and (2) $\lambda_{min}(\tilde{\Sigma}_1) \geq k_1 N^{1-2\alpha_1}$ and $\lambda_{min}(\tilde{\Sigma}_2) \geq k_1 N^{1-2\alpha_2}$, otherwise, provided that the function $f$ is bounded away from zero. Let $P_{h_L}$ be the probability law of the process $\pmb{h}_i$ under the hypothesis $h_L$ defined in \fr{f_j}. We consider two cases; case when $\nu_i=0$ and case when $\sigma_1=0$.  \\
{\bf Case when $\nu_i=0$ and $\alpha_2 \leq 1/2$}. Then, by Assumptions ${\bf A.2}$ and ${\bf A.3}$, the Kullback divergence can be written as 
\beqn
K(P_{h_{L}}, P_{\tilde{h}_{L}})&\leq& \frac{N\rho^2_L}{2\sigma^2_1}\lambda_{\max}\left[\left(\tilde{\Sigma}_2\right)^{-1}\right]\sum^{L-1}_{l=0}|\tilde{\omega}_l-\omega_l|^2{\int^{b}_0\varphi^2_l(t)g(t)dt}\nonumber\\
&\leq& \frac{N\rho^2_L}{2\sigma^2_1}\left[\left(\lambda_{\min}[\tilde{\Sigma}_2]\right)^{-1}\right]\sum^{L-1}_{l=0}|\tilde{\omega}_l-\omega_l|^2{\int^{b}_0\varphi^2_l(t)g(t)dt}\nonumber\\
& \leq & \frac{N\rho^2_L}{k_2\sigma^2_1N^{1-2\alpha_2}} 2L \max_{l \leq L-1}|\varphi_l(t)|^2=CN^{2\alpha_2}A^2L^{-(2s+1)}L.
\eeqn
Now, to apply Lemma \ref{lem:Bunea}, choose 
\be \label{klblam}
A^2 L^{-(2s+1)}LN^{2\alpha_2}\leq \pi_0L,
\ee
which gives
\be \label{Lop}
L = C\left[{A^2N^{2\alpha_2}}\right]^{\frac{1}{2s +1}}.
\ee
{\bf Case when $\sigma_1=0$ and $\alpha_1 \leq 1/2$}. Then, the Kullback divergence gives
\beqn
K(P_{h_{L}}, P_{\tilde{h}_{L}})&\leq& \frac{N\rho^2_L}{2}\lambda_{\max}\left[\left(\tilde{\Sigma}_1\right)^{-1}\right]\sum^{L-1}_{l=0}|\tilde{\omega}_l-\omega_l|^2{\int^{b}_0\varphi^2_l(t)g(t)dt}\nonumber\\
&\leq& \frac{N\rho^2_L}{2}\left[\left(\lambda_{\min}[\tilde{\Sigma}_1]\right)^{-1}\right]\sum^{L-1}_{l=0}|\tilde{\omega}_l-\omega_l|^2{\int^{b}_0\varphi^2_l(t)g(t)dt}\nonumber\\
& \leq & \frac{N\rho^2_L}{k_1N^{1-2\alpha_1}} 2L \max_{l \leq L-1}|\varphi_l(t)|^2=CN^{2\alpha_2}A^2L^{-(2s+1)}L.
\eeqn
Now, to apply {\bf Lemma \ref{lem:Bunea}}, choose 
\be \label{klblam}
A^2 L^{-(2s+1)}LN^{2\alpha_1}\leq \pi_0L,
\ee
which gives  
\be \label{Lop}
L = C\left[{A^2N^{2\alpha_1}}\right]^{\frac{1}{2s +1}}. 
\ee
Notice that cases {\bf Case when $z_i=0$ and $\alpha_1 > 1/2$} and {\bf Case when $\nu_i=0$ and $\alpha_2 > 1/2$} the matrices $\tilde{\Sigma}_1$ and $\tilde{\Sigma}_2$ have finite eigenvalues (they do not depend on $N$) can be dealt with in a similar fashion as to the i.i.d case, so we skip them. 
To complete the proof, keep in mind that the highest of the lower bounds corresponds to 
\be
\tilde{L}=\min\left\{ \left[{A^2N^{2\alpha_1}}\right]^{\frac{1}{2s +1}}, \left[{A^2N^{2\alpha_2}}\right]^{\frac{1}{2s +1}}, \left[{A^2N}\right]^{\frac{1}{2s +1}} \right\}.  \ \Box
\ee
{\bf Proof of Lemma \ref{lem:Var}.} Notice that with $\widehat{\theta}_l$  defined in \fr{theta-h}, one has 
\be \label{alef}
\widehat{\theta}_l-\theta_l= \frac{1}{N}\sum^N_{i=1}\left[\eta_i\II\left(\Omega_l(i) \right)-\EE\left[\eta_i\II\left(\Omega_l(i) \right)\right]\right]- \EE\left[\eta_i\II\left(\Omega^c_l(i) \right)\right],
\ee
where the quantities $\eta_i=\left[y_i^2\frac{\varphi_l(t_i)}{g(t_i)}-\sigma^2\int^{b}_0\varphi_l(t)dt\right]$. Define the quantities $\Delta_i=[\eta_i\II\left(\Omega_l(i) \right)-\EE\left[\eta_i\II\left(\Omega_l(i) \right)\right]$ and notice they  are independent zero-mean random variables with variance 
\be \label{deli}
\EE\left[\Delta_i\right]^2\leq 2\EE\left[ y^4_i\frac{\varphi^2_l(t_i)}{g^2(t_i)}\right]+8\sigma^4=2\left[\EE[\varepsilon_i^4]\int f^4(t)\frac{\varphi^2(t)}{g(t)}dt +6\sigma^2\int f^2(t)\frac{\varphi^2(t)}{g(t)}dt+\sigma^4\EE[z^4_i]\int \frac{\varphi^2(t)}{g(t)}dt\right] +8\sigma^4=\sigma^2_0.
\ee
In addition, by Cauchy-Schwarz inequality and the Gaussian tail probability inequality, we want to  show that 
\be \label{bia}
\left(\EE\left[\eta_i\II\left(\Omega^c_l(i) \right)\right] \right)^2\leq \EE\left[\eta^2_i\right]\Pr\left(\Omega^c_l(i) \right)=o\left(N^{-1}\right). 
\ee
Bear in mind that conditional on the distribution $g$, the quantities $y_i=\varepsilon_i f(t_i)+\sigma z_i$ are $N(0, f^2(t_i) +\sigma^2)$. Therefore, by Assumptions ${\bf A.1}$ and ${\bf A.2}$ and equation $(2.5)$ of Muckenhoupf~(1970), we obtain
\beqns
\Pr\left(\Omega^c_l(i) \right)&\leq& \Pr\left(\frac{|\varphi_k(t_i)|y^2_i}{m_1(f^2(t_i) +\sigma^2)}\geq \frac{\sigma^2\int^{b}_0\varphi_l(t)dt\ + \sqrt{N/\ln(N)}}{f^2(t_i) + \sigma^2}\right)\\
& \leq & \Pr\left(\frac{C_0y^2_i}{m_1(f^2(t_i) +\sigma^2)}\geq \frac{\sigma^2\int^{b}_0\varphi_l(t)dt\ + \sqrt{N/\ln(N)}}{f^2(t_i) + \sigma^2}\right)\\
&=&\Pr\left(\frac{y_i}{\sqrt{(f^2(t_i) +\sigma^2)}}\geq \sqrt{\frac{m_1}{C_0}\frac{\sigma^2\int^{b}_0\varphi_l(t)dt\ + \sqrt{N/\ln(N)}}{f^2(t_i) + \sigma^2}}\right)\\
&\leq & C_2 \exp \left\{-\frac{m_1}{2C_0}\frac{\sqrt{N/\ln(N)}}{M_2^2+\sigma^2}\right\}.
\eeqns
To complete the proof, apply the expectation to the square of \fr{alef} and use results \fr{deli} and \fr{bia}. $\Box$\\
{\bf Proof of Lemma \ref{lem:Lardev}.}  In order to prove \fr{Largdev}, we make use of Bernstein inequality. 
    \begin{lemma}(Bernstein Inequality). \label{lem:bernineq}
Let $Y_i$, $i=1, 2, \cdots, N$,  be independent and identically distributed random variables with mean zero and finite variance $\sigma^2$, with $\|Y_i\| \leq \|Y\|_{\infty} < \infty$. Then,  
\be \label{prob-b}
\Pr\left(\left|N^{-1}\sum^N_{i=1}Y_i\right|>z\right) \leq 2 \exp\left\{-\frac{Nz^2}{2(\sigma^2+\|Y\|_{\infty}z/3)}\right\}.
\ee
\end{lemma}
Recall the notation used in the proof of {\bf Lemma \ref{lem:Var}}. Thus, since $\Var(\Delta_i) \leq \sigma^2_0<\infty$, and $\left|\Delta_i\right| \leq c_o\sqrt{N/\ln(N)}$, for $\gamma\rho>1$, Bernstein inequality gives 
\beqns
\Pr \left(| \widehat{\theta}_{l}-\theta_l |> \rho\lambda_l\right)\leq \Pr \left(\left| N^{-1}\sum^N_{i=1}\Delta_i \right|+\EE\left[\eta_i\II\left(\Omega^c_l(i) \right)\right]> \rho\lambda_l\right)&\leq& \Pr \left(\left| N^{-1}\sum^N_{i=1}\Delta_i \right|> (\gamma\rho-1)\sqrt{\ln(N)/N}\right)\\
&\leq& 2\exp\left\{ -\frac{(\gamma\rho-1)^2\ln(N)}{2(\sigma^2_0+c_o(\gamma\rho-1)/3)}\right\}\\
&=& 2 N^{-\frac{(\gamma\rho-1)^2}{2(\sigma_0^2+c_o/3(\gamma\rho-1))}}. \Box
\eeqns
{\bf Proof of Theorem \ref{th:upperbds-2}}.   Denote
\be \label{chijj}
\chi_{N}=\frac{\ln(N)}{A^2N}, \ \ M_o=[\chi_{N}]^{-\frac{1}{2s + 1}},
\ee
 and note that with the choice of $M$ and $\lambda_l$ given by \fr{Lev:Jad}  and \fr{thresh1}, respectively, the estimation error can be decomposed as $\mathbb{E} \| \widehat{h}_M-h \|^2
\leq\mathbb{E}_1 +\mathbb{E}_2 +\mathbb{E}_3+\mathbb{E}_4+\mathbb{E}_5$, where
 \beqn  
\mathbb{E}_1&=&  \sum^{M-1}_{l=1}\EE \left[\left| \widehat{ \theta}_{l}-\theta_{l}\right|^2 \II \left(   \left| \widehat{ \theta}_{l}-\theta_{l}  \right| > \frac{1}{2} \lambda_{l}\right)  \right], \label{r21}\\
\mathbb{E}_2&=& \sum^{M-1}_{l=1}\EE \left[\left| \widehat{ \theta}_{l}-\theta_{l}\right|^2 \II \left(  \left|  \theta_{l} \right| >  \frac{1}{2}\lambda_{l} \right)  \right],\label{r22}\\
\mathbb{E}_3&=&\sum^{M-1}_{l=1}{ \theta}_{l}^2 \Pr \left( \left| \widehat{ \theta}_{l}-\theta_{l}  \right| > \frac{1}{2} \lambda_{l} \right),\label{r31}\\
\mathbb{E}_4&=& \sum^{M-1}_{l=1}  { \theta}_{l}^2\II \left(  \left|  \theta_{l} \right| <  \frac{3}{2}\lambda_{l} \right).\label{r32}\\
\mathbb{E}_5&=& \sum^{\infty}_{l=M} { \theta}_{l}^2.\label{r33}
\eeqn
Then, by  \fr{eq11} and \fr{Lev:Jad}, \fr{r33} becomes 
\beqn \label{ee14}
\mathbb{E}_5
&=&O\left( \sum^{\infty}_{l=M} A^2 \left(l\vee 1\right)^{-2s} \right)=O\left(A^2 \left[\chi_{n}\right]^{{2s}}\right)=O\left(A^2 \left[\chi_{n}\right]^{\frac{2s}{2s+1}}\right).
\eeqn
Now, combining $\mathbb{E}_1$ and $\mathbb{E}_3$, and applying Cauchy-Schwarz inequality, the moments property of the Gaussian, Lemma \ref{lem:Var} with the choice $\tau > 2$, \fr{eq11} and \fr{Lev:Jad}, yields 
\be
\mathbb{E}_1 +\mathbb{E}_3 = O  \left( \frac{M}{N}\left[ \frac{1}{N}\right]^{{\tau}/{2}}+ A^2M\left[ \frac{1}{N}\right]^{{\tau}} \right)=  O \left(  \frac{1}{N} \right). \label{r21r31}
\ee
Now, combining $\mathbb{E}_2$ and $\mathbb{E}_4$ and using condition \fr{eq11} yields
\be  \label{r22r32}
 \Delta=\mathbb{E}_2 +\mathbb{E}_4=  O \left( \sum^{M-1}_{l=1}\min \left\{  { \theta}_{l} ^2,   A^2 \left[ \chi_{n}\right]\right\} \right).
\ee
Finally, $\Delta$ can be decomposed into the following components
\beqn  
 \Delta_1&=&  O \left( \sum^{M-1}_{l=M_0}A^2 \left(l\vee 1\right)^{-2s} \right)= O\left(A^2 \left[\chi_{n}\right]^{\frac{2s}{2s+1}}\right), \label{del1}\\
 \Delta_2&=&O \left( \sum^{M_0-1}_{l=1}A^2  \left[ \chi_{n}\right]\right)=O \left( M_0A^2  \left[ \chi_{n}\right]\right) =O\left(A^2 \left[\chi_{n}\right]^{\frac{2s}{2s+1}}\right). \label{del2}
 \label{del3}
\eeqn
Hence, combining \fr{ee14}, \fr{r21r31}, \fr{del1} and \fr{del2} completes the proof. $\Box$\\
{\bf Proof of Lemma \ref{lem:Var-lm}.} Notice that with $\widehat{\theta}_l$  defined in \fr{theta-h-lm}, and using property \fr{auto-cov}, the properties $\EE[z^4_i]=3$ and $\EE[z^2_iz^2_j]=\EE^2[z^2_i]+2\EE^2[z_iz_j]$, and the assumptions ${\bf A.1}$, ${\bf A.2}$ for $m_1 \geq 1$, as $N \rightarrow \infty$, one has
\beqn  
\EE\left[\widehat{\theta}_l-\theta_l\right]^2&=& \EE\left[ \frac{1}{N}\sum^N_{i=1}\left(\varepsilon^2_i f^2(t_i)+2\sigma \varepsilon_iz_i f(t_i)+\sigma^2z^2_i\right)\frac{\varphi_l(t_i)}{g(t_i)}\right]^2-\sigma^4\left(\int^b_0\varphi_l(t)dt\right)^2-\theta^2_l-2\sigma^2\theta_l\int^b_0\varphi_l(t)dt\nonumber\\
&=&\frac{3}{N}\left[\int^b_0 f^4(t)\varphi^2(t)/g(t)dt +\sigma^2\int^b_0 f^2(t)\varphi^2(t)/g(t)dt+\sigma^4\int^b_0\varphi^2_l(t)/g(t)dt\right]\nonumber\\
&+& \frac{1}{N^2}\sum_{i\neq j}\left[\theta_l^2\EE[\varepsilon^2_i \varepsilon^2_j]+2\sigma^2\theta_l\int^b_0 \varphi_l(t) dt +4c_1c_2\sigma^2|i-j|^{-\alpha_1-\alpha_2}\left(\int^b_0 f(t)\varphi_l(t)dt\right)^2\right]\nonumber\\
&+&\frac{1}{N^2}\sum_{i\neq j}\left[\sigma^4\EE[z^2_iz^2_j]\left(\int^b_0 \varphi_l(t) dt\right)^2 \right]- 2\sigma^2\theta_l\int^b_0 \varphi_l(t)dt-\sigma^4\left(\int^b_0\varphi_l(t)dt\right)^2-\theta^2_l \nonumber\\
&\leq& \frac{3}{N}\left[M^2_2+\sigma^2\right]^2+ \frac{2}{N^2}\sum_{i\neq j}\left[M^2_2c_1|i-j|^{-\alpha_1} +\sigma^2c_2|i-j|^{-\alpha_2} \right]^2.
\eeqn
Notice that in the last line, if both $\alpha_1$ and $\alpha_2$ are greater than 1/2, the first term with dominate, otherwise, the variance will be bounded by the larger of ${N^{-2\alpha_j}}$, $j=1, 2$. $\Box$\\
{\bf Proof of Lemma \ref{lem:Lardev-lm}.} Below, we use a combination of Lemma 2 in \cite{ben2}, which is an adaptation of Hanson-Wright inequality to matrices, and large deviation result that was developed in \cite{com1} and further improved in \cite{gen}  which states that for any $x>0$, if $\xi_n$ is a zero mean Gaussian vector with independent elements, and $Q$ is nonnegative definite matrix, then 
\be \label{QLarge-D}
\Pr\left(\xi_n^TQ\xi_n> \sigma^2\left[\sqrt{(tr(Q))}+\sqrt{x\rho^2_{\max}(Q)}\right]^2\right)\leq e^{-x}.
\ee
Let $F$ and $\Phi$ be the $N$-dimensional diagonal matrices whose diagonal elements are $f(t_1), f(t_2), \cdots, f(t_N)$, and $\varphi_l(t_1)/g(t_1), \varphi_l(t_2)/g(t_2), \cdots, \varphi_l(t_N)/g(t_N)$, respectively. Then, 
\beqn
\Pr \left(| \widehat{\theta}_{l}-\theta_l |> \rho\lambda_l\right)&\leq& \Pr \left(\left|2\pmb{ \varepsilon}^T_NF\Phi F\pmb{ \varepsilon}_N-N\theta_l\right| > \frac{\rho\gamma}{2}\frac{\ln{(N)}}{N^{\alpha_1-1}}\right)+\Pr \left(\left|2\pmb{z}^T_N\Phi \pmb{z}_N-N\int\varphi_l(t)dt\right| >\frac{\rho\gamma}{2\sigma^2}\frac{\ln{(N)}}{N^{\alpha_2-1}}\right)\nonumber\\
&=& P_1+P_2.
\eeqn
For the first term and if $\theta_l >0$, we apply \fr{QLarge-D}. Therefore, with \fr{epseta} and \fr{LRD}, as $N \rightarrow \infty$, one has 
\be
\mathrm{Tr}\left(A_1^TF\Phi FA_1\right)=\mathrm{Tr}\left(F\Phi FA_1A_1^T\right)=\mathrm{Tr}\left(F\Phi F\Sigma_1\right)=\sum^N_{i=1}f^2(t_i)\varphi_l(t_i)/g(t_i)\leq \frac{N}{bm_1}\theta_l,
\ee
and
\be
\rho^2_{\max}\left(A_1^TF\Phi FA_1\right)=\rho^2_{\max}\left(F\Phi F\Sigma_1\right)\leq \rho^2_{\max}\left(F\Phi F\right)\rho^2_{\max}\left(\Sigma_1\right)\leq \frac{M^2_2}{m_1}\max_{i}|\varphi_l(t_i)|c^{(1)}_2N^{1-\alpha_1}=\frac{M^2_2}{m_1}\pi_0c^{(1)}_2N^{1-\alpha_1}\ee
Therefore, if $bm_1 \geq 4\Var(\eta^{(1)}_i)$, 
\beqn
\Pr \left(\left|2\pmb{ \varepsilon}^T_NF\Phi F\pmb{ \varepsilon}_N-N\theta_l\right| > \frac{\rho\gamma}{2}\frac{\ln{(N)}}{N^{\alpha_1-1}}\right)&\leq& \Pr\left(\pmb{ \varepsilon}^T_NF\Phi F\pmb{ \varepsilon}_N>\frac{bm_1}{4}\left[\sqrt{\frac{N\theta_l}{bm_1}}+\sqrt{ \frac{\rho\gamma}{2bm_1}\frac{\ln{(N)}}{N^{\alpha_1-1}}}\right]^2\right).\nonumber
\eeqn
\fr{QLarge-D} is applied then by taking $x=\frac{\rho\gamma\ln{(N)}}{c^{(1)}_2\pi_0bM^2_2}$. Now, if $\theta_l <0$, then we apply Hanson-Wright inequality from \cite{RuV} to
\be
\Pr \left(\left|2\pmb{ \varepsilon}^T_NF\Phi F\pmb{ \varepsilon}_N-N\theta_l\right| > \frac{\rho\gamma}{2}\frac{\ln{(N)}}{N^{\alpha_1-1}}\right) \leq \Pr \left(\left|\pmb{ \varepsilon}^T_NF\Phi F\pmb{ \varepsilon}_N-N\theta_l\right| > \frac{\rho\gamma}{4}\frac{\ln{(N)}}{N^{\alpha_1-1}}\right),\nonumber
\ee
with matrix $B=A_1^TF\Phi FA_1$ having Frobenius norm 
\be
\|A_1^TF\Phi FA_1\|^2_F\leq \mathrm{Tr}^2\left(A_1^TF\Phi FA_1\right)= \mathrm{Tr}^2\left(F\Phi F\Sigma_1\right)\leq N\frac{M^4_2}{b^2m^2_1},\ \ as\ N \rightarrow \infty. \nonumber
\ee
Hence, applying Hanson-Wright inequality yields 
\be
\Pr \left(\left|\pmb{ \varepsilon}^T_NF\Phi F\pmb{ \varepsilon}_N-N\theta_l\right| > \frac{\rho\gamma}{4}\frac{\ln{(N)}}{N^{\alpha_1-1}}\right) \leq 2 \exp\left\{-\frac{c_om_1\gamma_1\rho\ln{(N)}}{4k^2_1\pi_0c^{(1)}_2}\right\}.
\ee
In a similar fashion, one can evaluate $P_2$ taking into consideration whether $\int^b_0\varphi_l(t)dt$ is positive, in which case we use \fr{QLarge-D}, or negative in which case we apply  Hanson-Wright inequality. $\Box$\\
{\bf Proof of Theorem \ref{th:upperbds-lm}}. The proof is similar to that of {\bf Theorem \ref{th:upperbds-2}}, so we skip it.  $\Box$


\begin{thebibliography}{99}

 \bibitem{ben2} Benhaddou, R. (2018a), 'Laplace deconvolution  with dependent errors : a Minimax Study', {\it Journal of Nonparametric Statistics}, {\bf 30(4)}, 1032-1048.
 
  \bibitem{ben6} Benhaddou, R. (2018b), 'Minimax lower bounds for the simultaneous wavelet deconvolution  with fractional Gaussian noise and unknown kernels', {\it Statistics and Probability Letters}, {\bf 140}, 91-95.
  
   \bibitem{ben4} Benhaddou, R. (2016), 'Deconvolution  model with fractional Gaussian noise: a Minimax Study', {\it Statistics and Probability Letters}, {\bf 117}, 201-208.
   
    \bibitem{ben1} Benhaddou, R., Kulik, R., Pensky, M., Sapatinas, T. (2014), 'Multichannel deconvolution with long-range dependence: a minimax Study', {\it  Journal of Statistical Planning and Inference}, {\bf 148}, 1-19.

  \bibitem{ben5} Benhaddou, R., Liu, Q. (2019), 'Anisotropic functional deconvolution with long-memory noise: the case of a multi-parameter fractional Wiener sheet', {\it Journal of Nonparametric Statistics}, {\bf 31(3)}, 567-595.

 \bibitem{ben3} Benhaddou, R., Pensky, M., Rajapakshage, R. (2019), 'Anisotropic functional Laplace deconvolution', {\it Journal of Statistical Planning and Inference}, {\bf 199}, 271-285.
 
   \bibitem{bong} { Bongioanni, B., Torrea, J. L. (2009), 'What is a Sobolev space for the Laguerre function systems ?', {\it Studia Mathematica}, {\bf 192(2)}, 147-172}.
 
\bibitem{br} Brown, L.D., Levine, M. (2007), 'Variance estimation in nonparametric regression via the difference sequence method', {\it The Annals of Statistics}, {\bf 35(5)}, 2219-2232.  

    \bibitem{bun}
Bunea, F., Tsybakov, A. \& Wegkamp, M.H. (2007),
'Aggregation for Gaussian regression',
{\it Annals of Statistics,} {\bf 35}, 1674--1697.

  \bibitem{cai} Cai, T.T., Wang, L. (2008), 'Adaptive  variance  function  estimation  in  heteroscedastic nonparametric regression', {\it The Annals of Statistics}, {\bf 36(5)}, 2025-2054  
  
    \bibitem{ches} Chesneau, C., El Kolei, S., Kou, J., Navarro, F. (2020), 'Nonparametric estimation in a regression model with additive and multiplicative noise', {\it Journal of Computational and Applied Mathematics}, {\bf 380}.    
    
\bibitem{chichi} Chichignoud, M. (2012), 'Minimax and minimax adaptive estimation in multiplicative regression: locally Bayesian approach', {\it Probability Theory and Related Fields}, {\bf 153(3-4)}, 543-586.
  
   \bibitem{com1}
Comte, F. (2001),
'Adaptive Estimation of the Spectrum of a Stationary Gaussian Sequence',
{\it Bernoulli}, {\bf 7},   267--298.

\bibitem{com5}
Comte, F.,    Cuenod, C.-A.,    Pensky, M.,   Rozenholc, Y. (2017),
'Laplace deconvolution on the basis of  time domain data  
and its application to Dynamic Contrast Enhanced Imaging',
{\it Journal of the Royal Statistical Society, Ser.B}, {\bf 79},   69--94.

\bibitem{com3}
Comte, F., Dedecker, J., Y., Taupin, M.L. (2008),
'Adaptive density deconvolution with dependent inputs',
{\it Mathematical Methods in Statistics}, {\bf 17},   87--112.

\bibitem{com4}
Comte, F., Genon-Catalot, V. (2015),
'Adaptive Laguerre density estimation for mixed Poisson models',
{\it Electronic Journal of Statistics}, {\bf 9(1)},   1113--1149.


 \bibitem{gen} Gendre, X. (2014), 'Model Selection and Estimation of a Component in Additive  Regression', {\it ESAIM: Probability and Statistics}, {\bf 18}, 77-116.

 \bibitem{HT} Hardle, W., and Tsybakov, A. (1997), 'Local polynomial estimators of  volatility function in nonparametric regression', {\it The Journal of Econometrics}, {\bf 81(1)}, 223-242.
 
  \bibitem{Hu} { Huang, P., Pi, Y., and Progri, I. (2013), 'GPS Signal Detection under Multiplicative and Additive Noise', {\it The Journal of Navigation}, {\bf 66}, 479-500}.
 
   \bibitem{kl1} Kulik, R., Raimondo, M. (2009), 'Wavelet regression in random design with heteroskedastic dependent errors', {\it  The Annals of Statistics}, {\bf 37(6)}, 3396-3430.
    
   \bibitem{kl} Kulik, R., Wichelhaus, C. (2011), 'Nonparametric conditional variance and error density estimation in regression models  with dependent errors and predictors', {\it  Electronic Journal of Statistics}, {\bf 5}, 856-898.
   
\bibitem{mab}
Mabon, G. (2016),
'Adaptive deconvolution of linear functionals on nonnegative real line', 
{\it Journal of Statistical Planning and  Inference}, {\bf 178},  1-23.
 
 \bibitem{abram} Muckenhoupt, B. (1970), 'Mean convergence of Hermite and Laguerre series II',
{\it Translations of the American  Mathematical Society}, {\bf 147}, 433-460.

\bibitem{RuV} Rudelson, M., Vershynin, R. (2013), 'Hanson-Wright inequality and sub-Gaussian concentration', {\it Electronic Communications in Probability}, {\bf 18(82)}, 1-19.

\bibitem{tsybakov}
Tsybakov, A.B. (2008),
{\it Introduction to Nonparametric Estimation}, Springer, New York.

\bibitem{vareschi}
Vareschi, T. (2015),
'Noisy Laplace deconvolution with error in the operator', 
{\it Journal of Statistical Planning and  Inference}, {\bf 157-158},  16-35.

\bibitem{wan1} Wang, Y. (1996), 'Functional estimation via wavelets shrinkage for Long-memory Data', {\it Annals of Statistics}, {\bf 24}, 466-484.

\bibitem{wan2} Wang, Y. (1997), 'Minimax Estimation via Wavelets for Indirect Long-memory Data', {\it Journal of Statistical Planning and  Inference}, {\bf 1}, 45-55.

\bibitem{wish} Wishart, J. M. (2013), 'Wavelet deconvolution in a periodic setting with long-range dependent  errors', {\it Journal of Statistical Planning and  Inference}, {\bf 5}, 867-881.


\end{thebibliography}
\end{document}